%% file: main.tex
\documentclass{elsarticle}
\usepackage{geometry}
\input{preamble}
\title{Matrix-Free Methods for Finite-Strain Elasticity: Automatic Code Generation with No Performance Overhead}
\author[IWR]{Micha\l{} Wichrowski\corref{cor1}}
\author[IPPT]{ Mohsen Rezaee-Hajidehi}
\author[UL]{ Jože Korelc}
\author[RUB]{ Martin Kronbichler}
\author[IPPT]{ Stanis\l{}aw Stupkiewicz }

\cortext[cor1]{Corresponding author, mt.wichrowsk@uw.edu.pl}

\address[IWR]{Interdisciplinary Center for Scientific Computing, Heidelberg University, Heidelberg, Germany}
\address[IPPT]{Institute of Fundamental Technological Research, Polish Academy of Sciences, Warsaw, Poland}
\address[UL]{University of Ljubljana, Faculty of Civil and Geodetic Engineering, Slovenia}
\address[RUB]{Applied Numerics, Faculty of Mathematics, Ruhr University Bochum, Germany}
\date{May, 2025}

\usepackage{color}
\usepackage{bbm}

\newcommand{\plotscale}{.5}

\begin{document}

\begin{abstract}
	This study explores matrix-free tangent evaluations in finite-strain elasticity with the use of
	automatically-generated code for the quadrature-point level calculations. The code generation is done via
	automatic-differentiation (AD) with AceGen. We compare hand-written   and AD-generated codes under two
	computing
	strategies: on-the-fly evaluation and caching intermediate results. The comparison reveals that the
	AD-generated code
	achieves superior performance in matrix-free computations.
\end{abstract}

\maketitle

\section{Introduction}
\label{sec:intro}
Matrix-free methods	 rely on the elegance of well-optimized loops to evaluate the action of a linear operator on a
vector without explicitly storing matrix entries~\cite{orszag1980spectral, brown2010efficient, cantwell2011from,
	kronbichler2012generic,
	gmeiner2016quantitative, kronbichler2019fast, abdelfattah2021gpu, kolev2021efficient}.
While these methods significantly boost computational performance, including applications in
solid mechanics~\cite{davydov2020matrix, brown2022performance, schussnig2024matrix}, their adoption is hindered by
the formidable challenge of deriving and implementing complex tangent operators. This study
explores the applicability of automatically generated codes within matrix-free non-linear solvers, focusing on
finite-strain elasticity.

Typical finite-element computations involve operating on large sparse matrices and are, as a consequence, memory-bound.
This implies that the throughput (time per unknown) of sparse matrix--vector products is limited by main memory access
speed,
leaving a vast majority of computational resources fallow.
In~\cite{brown2022performance,williams2009roofline}
it was estimated that for a modern CPU, iterative sparse linear solvers saturate memory bandwidth at less than 2\% of
the processor's  theoretical arithmetic computing power. As computing capabilities continue to outgrow memory
bandwidth~\cite{gholami2024ai}, this gap is further exacerbated.
At the same time, higher-order elements, with denser global tangent matrix, can decrease throughput in
classical matrix-based solvers.
Consequently, linear
elements are often preferred for large-scale simulations, despite their lower accuracy per degree of
freedom~\cite{brown2022performance,	schneider2022large, duster2003p} and susceptibility to
locking~\cite{wriggers2008nonlinear}.

Matrix-free operator evaluation alleviates the
memory bottleneck by avoiding the explicit storage and manipulation of sparse matrices~\cite{kolev2021efficient,
	may2014ptatin3d, kronbichler2012generic}.
The most general matrix-free approach involves computing finite-element integrals in weak forms on the fly.
Leveraging the tensor-product structure of the finite-element basis functions through sum
factorization~\cite{melenk2001fully} allows to reduce the arithmetic load of computing the integrals for every operator
evaluation.
Additional optimizations can further reduce the computational effort associated with the one-dimensional
interpolations in sum factorization, e.g., by the even--odd decomposition proposed by~\cite{solomonoff1992fast} and
analyzed in~\cite{kronbichler2019fast}. Due to the higher arithmetic intensity, hardware-specific optimizations such as
the use of vectorization (SIMD
instructions) have been developed.	Some optimization strategies specific to triangles, tetrahedra, and prismatic
elements have been also proposed~\cite{moxey2020efficient}.  A review of high-performance solvers can be found
in~\cite{fischer2020scalability},
which explores performance and space-time trade-offs for compute-intensive kernels of large-scale numerical solvers for
PDEs. For solving the linear systems arising from elliptic partial differential equations, matrix-free methods are
usually combined with multigrid
solvers~\cite{kronbichler2018performance,Clevenger2020Flexible}, leading to overall linear complexity with respect to
the number of degrees of freedom.

Employing the matrix-free method for solving problems of finite-deformation solid mechanics poses significant
challenges. A primary challenge is the complexities involved in the derivation of the tangent operator, as the
proposed quadrature-based approach needs to evaluate the constitutive terms in every evaluation step. A successful
implementation of the matrix-free method for a (compressible) neo-Hookean hyperelastic model was presented
by Davydov et al.~\cite{davydov2020matrix}. They proposed different caching strategies to evaluate the resulting
constitutive relations, using scalar quantities, second-order tensors or a
fourth-order tensor, with the aim to optimize the performance of the matrix-free method.
Nevertheless, their implementation relied on the explicit evaluation of the tangent operator, a requirement that
limits its applicability to more complex models.
More recently, in~\cite{wichrowski2023exploiting}, a transient
fluid--structure interaction solver  was presented that involves an incompressible Mooney--Rivlin solid; however, the
computations were performed without the need for the tangent operator.	The contribution~\cite{brown2022performance}
proposed GPU algorithms for structural mechanics, whereas~\cite{schussnig2024matrix} presented a matrix-free
finite element solver for finite-strain elasticity with an \textit{hp}-multigrid preconditioner. The latter focused on
storage
strategies
to balance compute load against memory
access for fiber-reinforced models. Both studies demonstrated significant speed-ups for higher polynomial degrees in a
complex geometry.

The present work employs the general-purpose infrastructure provided by the
\texttt{deal.II} finite element library~\cite{dealii2019design,dealII95}, similar to the previous
contributions~\cite{davydov2020matrix, schussnig2024matrix}. Here,
the operations at quadrature points for both the residual and the associated tangent operator are specified
in the library code. Conversely, the interpolation of values or gradients to
quadrature points, summation for quadrature, the loop over mesh elements and
the exchange of data between different processes in a parallel computation
are using library code. This enables the use of code tuning and data access
optimizations from previous contributions.

Automatic differentiation (AD) offers a powerful approach to generating efficient finite-element codes by reducing
manual intervention in differentiation and coding
\cite{korelc1997automatic,bischof2003extending,griewank2008evaluating}.
This minimizes human error in deriving the residual vector and
tangent operator, improving the reliability and efficiency of generated codes. Moreover, this automation can bring
significant time savings in the code development phase. In this study, we evaluate the performance of
matrix-free implementation of AD-generated finite-element codes for neo-Hookean hyperelasticity models using the
AceGen system (http://symech.fgg.uni-lj.si/), see also \cite{korelc2002multi,korelc2016}.
AceGen employs a hybrid symbolic--numerical approach to automate the finite element
method. It leverages the symbolic and algebraic capabilities of the general computer algebra system
Mathematica~\cite{Mathematica}, combined with automatic differentiation and code generation/optimization, to create
finite-element user subroutines.

AD tools can be broadly categorized by their implementation approach (operator overloading vs.\ source-to-source
transformation) and by the mode of applying the chain rule (forward vs.\ backward/inverse). Pure operator overloading
often suffers from low numerical efficiency. Consequently, most AD tools utilize some form of operator overloading to
generate code in a simplified intermediate language, followed by a source-to-source transformation.
AceGen uses a variant of the source-to-source implementation of AD in both backward and forward
modes, specifically optimized for generating finite-element subroutines.
A key feature is that the codes for a function and its derivatives are merged and optimized together. This means that
when calculating higher derivatives or applying directional derivatives repeatedly on the same function, only one
optimized subroutine is generated. These special features of AceGen's AD are ideal for an efficient automatic
generation
of subroutines for the matrix-free solution of linear systems. All operations with tensors and matrices are performed
and optimized component by component, completely avoiding function calls. In addition, AceGen's AD allows modifications
of
the chain rule through an AD exception mechanism~\cite{korelc2009automation}, which facilitates the
differentiation-based description of mechanical models.

%% SS: deleted old paragraph - this was a repetition

In~\cite{brown2022performance}, an attempt to use automatic differentiation through
Enzyme~\cite{moses2021reverse} for finite-strain elasticity yielded unsatisfactory results, with performance
over 30\% slower than hand-written code due to Enzyme's limitation to built-in types and incompatibility with
vectorized operations.

Our investigation focuses on a compressible neo-Hookean hyperelastic model, with efficient matrix-free implementations
available as a baseline~\cite{davydov2020matrix}. We compare our AD-based implementation with different caching-based
implementations from that work
and a conventional implementation using a sparse iterative solver~\cite{heroux2005overview}. Our
investigation focuses on
the trade-off between caching and computing, the overhead associated with AD, and the feasibility of storing partial
results while maintaining a general solid mechanics solver. Our results show that there is no overhead caused by AD,
and that the AceGen-generated code stands out as the best matrix-free implementation, at the same time being a general
approach, not limited to a specific constitutive model. We also evaluate our AD-generated code against the hand-crafted
implementation from \cite{schussnig2024matrix}, where another variant of hyperelastic model (with isochoric--volumetric
split of the energy) is considered, again demonstrating superior performance.

The remainder of this paper is organized as follows. We detail the formulation of the non-linear problem and the
solution procedure using the Newton method in Section~\ref{sec:nonlinear}. In particular, we describe the matrix-free
evaluation of the tangent operator and sum factorization in Section~\ref{sec:mf_workflow}, and the matrix-free
implementation using AD and caching strategies in Section~\ref{sec:pointwise}. We assess and compare the performance of
different matrix-free implementations in Section~\ref{sec:performance}. Finally, we wrap up the paper by outlining the
concluding remarks and potential future directions.

\section{The nonlinear problem}
\label{sec:nonlinear}

\subsection{Problem formulation}

We consider a hyperelastic body  that in the reference configuration occupies the domain $\Omega \subset
	\mathbb{R}^d$
with the boundary partitioned into a Dirichlet part $\Gamma_{\rm D}\subset\partial\Omega$ and a Neumann part
$\Gamma_{\rm N}\subset\partial\Omega$. We assume that the Dirichlet boundary $\Gamma_{\rm D}$ is fixed, while a
conservative
surface traction $\bfm{T}^\ast$ is applied on the Neumann boundary
$\Gamma_{\rm N}$. The deformation of the body is described  by the
mapping $\bm{\upvarphi}: \Omega \to \mathbb{R}^d$ that links the reference configuration $\Omega$ to the current
configuration $\omega$,
i.e., $\bm{\upvarphi}(\Omega )=\omega$.
The mapping $\bm{\upvarphi}$ is assumed to be a function having sufficient regularity for the weak formulation and
satisfying $J:=\det\bfm{F} > 0$ almost everywhere to preserve material orientation. The displacement field $\bfm{u}$ is
defined as the difference between the position in the deformed and reference configurations, i.e., $\bfm{u}(\bfm{X}) =
	\bm{\upvarphi}(\bfm{X}) - \bfm{X}$, where $\bfm{X}$ is the position vector in the reference configuration. The
deformation gradient
\begin{equation}
	\bfm{F} := \Grad \bm{\upvarphi}=\bfm{I} + \Grad \bfm{u}
\end{equation}
is the main kinematic quantity
in the finite-deformation setting. We use $\Grad(\cdot)$ to denote the gradient in the reference configuration.

The elastic response of the material is described by the strain energy density function $\Psi$, which is a
differentiable function of the deformation gradient $\bfm{F}$, or alternatively, of the right Cauchy--Green tensor
$\bfm{C}=\bfm{F}^{\rm T}\cdot\bfm{F}$. In this setting, we define the potential energy functional
\begin{equation}
	\mathcal{E}(\bfm{u}) := \int_\Omega \Psi(\bfm{F}) \; \Dd V
	- \int_{\Gamma_{\rm N}} \bfm{T}^\ast \cdot \bfm{u} \; \Dd S .
\end{equation}

We seek a displacement field $\bfm{u} \in \mathbb{V}$ that renders the first variation of $\mathcal{E}$ equal to zero
for all
admissible variations $\delta\bfm{u} \in \mathbb{V}$, where $\mathbb{V} = \{\bfm{v} \in H^1(\Omega)^d :
	\bfm{v}|_{\Gamma_{\rm D}} = \bfm{0}\}$ is the space of admissible displacements.
The corresponding stationarity condition yields the weak form, i.e., the virtual work principle,
\begin{align}
	\mathcal{F}(\bfm{u},\delta\bfm{u}) := \DD_{\delta\bfm{u}} \mathcal{E} & =
	\int_\Omega \frac{\partial\Psi}{\partial\; \bfm{\Grad u}} : \Grad\delta\bfm{u} \, \Dd V - \int_{\Gamma_{\rm N}}
	\bfm{T}^\ast
	\cdot
	\delta\bfm{u} \; \Dd S
	\nonumber
	\\
	\label{eq:F}
	                                                                      & = \int_\Omega \bfm{P} :
	\Grad\delta\bfm{u} \;
	\Dd V - \int_{\Gamma_{\rm N}} \bfm{T}^\ast \cdot \delta\bfm{u} \; \Dd S = 0 \qquad \forall \, \delta\bfm{u} ,
\end{align}
where $\DD_{\delta\bfm{u}} \mathcal{E}$ denotes the Gâteaux derivative in the direction $\delta\bfm{u}$, and
$$\bfm{P}:=\frac{\partial\Psi }{\partial\bfm{F}}= \frac{\partial\Psi}{\partial\; \bfm{\Grad u}}$$
is the first Piola--Kirchhoff stress tensor.

%%%%%%%%%%%%%%%%%%%%%%%%%%%%%%%%%%

To solve Eq.~\eqref{eq:F} for the unknown displacement $\bfm{u}$ we use Newton's method.
Given a current approximation of the displacement $\bar{\bfm{u}}$, we compute a Gateaux derivative of the functional
with respect to the displacement in direction $\Delta\bfm{u}$,
\begin{equation}
	\mathcal{K}(\bar{\bfm{u}};\Delta\bfm{u},\delta\bfm{u}) :=
	\DD_{\Delta\bfm{u}} \mathcal{F}(\bar{\bfm{u}},\delta\bfm{u}) .
\end{equation}
We compute the correction $\Delta\bfm{u}$ by solving the linearized problem
\begin{equation}
	\mathcal{F}(\bar{\bfm{u}}+\Delta\bfm{u},\delta\bfm{u}) \approx
	\mathcal{F}(\bar{\bfm{u}},\delta\bfm{u}) +
	\mathcal{K}(\bar{\bfm{u}};\Delta\bfm{u},\delta\bfm{u}) = 0 \qquad  \forall \, \delta\bfm{u} .
\end{equation}
The operator $\mathcal{K}(\bar{\bfm{u}};\Delta\bfm{u},\delta\bfm{u})$ is the tangent operator, bilinear with
respect to $\Delta\bfm{u}$ and $\delta\bfm{u}$. For conservative loading, it is given by
\begin{equation}
	\label{eq:DF}
	% \DD_{\Delta\bfm{u}} \mathcal{F}(\bar{\bfm{u}},\delta\bfm{u}) 
	\mathcal{K}(\bar{\bfm{u}};\Delta\bfm{u},\delta\bfm{u})
	= \int_\Omega \DD_{\Delta\bfm{u}}\bfm{P} : \Grad
	\delta\bfm{u} \; \Dd V = \int_\Omega \Grad\Delta\bfm{u} : \mathbb{L} : \Grad \delta\bfm{u} \; \Dd V ,
\end{equation}
where $\mathbb{L}$ is the fourth-order tangent stiffness tensor with the major symmetry ($L_{iAjB}=L_{jBiA}$),
\begin{equation}
	\label{eq:L}
	\mathbb{L} := \frac{\partial\bfm{P}}{\partial\bfm{F}} =
	\frac{\partial^2\Psi}{\partial\bfm{F}\otimes\partial\bfm{F}} .
\end{equation}

Both the weak form~\eqref{eq:F} and the tangent operator~\eqref{eq:DF} can be evaluated in the current configuration.
Specifically, the volume integral in Eq.~\eqref{eq:F} can be equivalently expressed as
\begin{equation}
	\label{eq:F:cur}
	\int_\Omega \bfm{P} : \Grad\delta\bfm{u} \; \Dd V = \int_\Omega \bm{\tau} : \grads \delta\bfm{u} \; \Dd V =
	\int_\omega \bm{\sigma} : \grads \delta\bfm{u} \; \Dd v ,
\end{equation}
where $\bm{\tau}=\bfm{P}\cdot\bfm{F}^{\rm T}$ is the Kirchhoff stress tensor, $\bm{\sigma}=\bm{\tau}/J$ is the Cauchy
stress tensor, and
$\Dd v=J\Dd V$. The gradient in the current configuration is denoted
by $\grad(\cdot)$, and $\grads\,(\cdot)$ denotes its symmetric part. Both $\bm{\tau}$ and $\bm{\sigma}$ are symmetric,
hence only the symmetric part of the gradient of $\delta\bfm{u}$ is involved.

Following~\cite{davydov2020matrix}, see also~\cite{wriggers2008nonlinear}, the tangent operator in Eq.~\eqref{eq:DF}
can be transformed to the current configuration as
\begin{equation}
	\label{eq:DF:cur}
	\mathcal{K}(\bar{\bfm{u}};\Delta\bfm{u},\delta\bfm{u})
	= \int_\omega \grads \Delta\bfm{u} : \bbc : \grads
	\delta\bfm{u} \, \Dd v + \int_\omega \grad \delta\bfm{u} : \big( \grads \Delta\bfm{u} \cdot \bm{\sigma} \big)
	\, \Dd v
	,
\end{equation}
where the first term is the material part of the tangent operator and the second term is the geometric part. Note that
the respective formula~(10) in~\cite{davydov2020matrix} is expressed through integrals over the reference configuration
$\Omega$, hence $J\bbc$ and $\bm{\tau}=J\bm{\sigma}$ are used instead of $\bbc$ and $\bm{\sigma}$, respectively.
The fourth-order spatial elasticity tensor $\bbc$ possesses both the minor and major symmetries
($c_{ijkl}=c_{jikl}=c_{klij}$) that we exploit in the implementation.
It is the push-forward of the material tangent stiffness tensor $\mathbb{C}$,
\begin{equation}
	\label{eq:Jc}
	J \bbc = \chi (\mathbb{C}) , \qquad
	\mathbb{C} := 4 \frac{\partial^2 \Psi(\bfm{C})}{\partial\bfm{C}\otimes\partial\bfm{C}} ,
\end{equation}
where $\Psi$ is now considered as a function of the right Cauchy--Green tensor $\bfm{C}:=\bfm{F}^{\rm T}\cdot\bfm{F}$,
and $\chi(\cdot)$
denotes the push-forward operation such that $\chi(\mathbb{C})_{ijkl}=F_{iA}F_{jB}F_{kC}F_{lD}C_{ABCD}$.

From the above formulation, it is evident that the crucial part of the evaluation involves
expressions with derivatives of the strain energy function $\Psi$, particularly its second derivatives as shown in
Eqs.~(\ref{eq:DF})--(\ref{eq:L}) and (\ref{eq:DF:cur})--(\ref{eq:Jc}).

\subsection{Numerical solution of the problem}

To solve this problem numerically, we apply the Finite Element Method (FEM) by introducing a mesh $\mesh_L$ that
subdivides the domain $\Omega$ into quadrilateral (2D) or hexahedral (3D) elements.   We define a finite-element space
$\mathbb{V}_L\subset \mathbb{V}$  of vector functions using the element $\mathbb{Q}_p$
of piecewise polynomials of degree  $p$ in each direction, see for instance~\cite{Ciarlet78}. Introducing the
finite-dimensional finite-element space reduces the problem to
finding the solution of a non-linear system of equations. This is achieved through the Newton method: at each
iteration, we solve a linear system of equations involving operator $\mathcal{K}(\bar{\bfm{u}}; \cdot,\cdot)$
to obtain a correction to the current solution.

We choose a basis  $\mathbb{B}=\{ \testFunc \in \mathbb{V}_L  \;	| \;  i=1,...,\dim \mathbb{V}_L\}$  using shape
functions $\testFunc$ in space $\mathbb{V}_L$. This choice  allows us to represent every element $\bfm{v}$ of space
$\mathbb{V}_L$ as a real vector $\fem{V}$, corresponding to the coefficients of the basis functions in the
expansion
of $\bfm{v}$. At each Newton iteration, for a given vector $\bar{\fem{U}}$ representing $\bar{\bfm{u}}$, the problem is
to find $\Delta\fem{U}$ with the corresponding $\Delta\bfm{u}$ such that
\begin{equation}
	\mathcal{K}(\bar{\bfm{u}}; \corr,\testFunc) = -\mathcal{F}(\bfm{u},\testFunc) \qquad \forall \testFunc
	\in \mathbb{B}.
\end{equation}

\subsection{Linear system and solver}
The linear system involving the tangent operator is typically large, and its solution often represents the most
time-consuming step of the procedure. For matrix-free approaches, the choice of the solver is restricted to
iterative methods, as the matrix is not explicitly formed. To achieve optimal convergence and performance, an effective
preconditioner is  crucial. In our case, we employ the conjugate gradient (CG) method in conjunction with a geometric
multigrid preconditioner, which is particularly well-suited for matrix-free implementations since it can be formulated
entirely in terms of matrix--vector products, simple operations such as the matrix diagonal, and grid transfer
operations.

%%Include this sentence nicely into the text below
Defining the multigrid iteration requires a hierarchy of problems as well as establishing level operators,
smoothers, and transfer operators. We build the multigrid procedure on the assumption that the finest mesh  $\mesh_L$
is a result of refining a coarse mesh so that nestedness is obtained:
\begin{gather}
	\mesh_0 \sqsubset \mesh_1 \sqsubset \dots \sqsubset \mesh_L.
\end{gather}
The symbol ``$\sqsubset$'' indicates that every cell of mesh $\mesh_{\ell+1}$ is obtained from a cell of mesh
$\mesh_{\ell}$ by refinement. On every level $\ell$, we define a finite-element space $\mathbb{V}_\ell$ in the same
manner
as on the
finest one. We assume the existence of transfer operators between these finite-element spaces, allowing us to
interpolate the current iterate $\bar{\bfm{u}}$ onto level $\ell$ to obtain $\bar{\bfm{u}}_\ell$, providing a
definition for tangent operator $\mathcal{K}(\bar{\bfm{u}}_\ell ; \cdot, \cdot)$ at each level, which will be used as a
level
operator within the multigrid procedure.

For the smoother, we utilize the inverse of the diagonal of the level tangents, using an iteration with
Chebyshev polynomials, as described in~\cite{adams2003parallel, kronbichler2012generic,kronbichler2018performance}.  As
a result,
the matrix--vector multiplication is the dominant operation in the smoother, which is in turn the dominant operation in
the
overall solver.
The optimizations enabled by the tensor-product structure
of $\mathbb{Q}_p$ elements and the low memory access of the matrix-free method are thus addressing the most
expensive step in the solver.

\subsection{Matrix-free evaluation of the tangent operator}
\label{sec:mf_workflow}

Applying the tangent operator in a matrix-free manner involves computing a result vector $\resultReal$ from an input
vector $\corrReal$ that represents the finite-element
%% function?
field $\corr$.
Each component of the result vector represents the application of the tangent operator to the input vector, tested with
a corresponding basis function, such that $\resultReal_i = \mathcal{K}(\bar{\bfm{u}}; \corr, \testFunc)$ for each basis
function $\testFunc \in \mathbb{B}$.
To compute these coefficients, we decompose the evaluation into cell-wise
contributions and apply numerical integration over each cell. For a cell $K \in \mathcal{T}$, the local contribution to
the $i$-th component is computed as:
\begin{equation}
	\resultReal_{K,i}=\int_K \Grad \corr :\mathbb{L} : \Grad {\bm{\phi}}_i \; \Dd V
	\approx \sum_q \Grad \corr :\mathbb{L} : \Grad\bm{\phi}_i \, J_q \, w_q.
\end{equation}
This numerical integration uses $(p+1)^d$ quadrature points, where $J_q$ denotes the Jacobian
determinant at point $q$, and $w_q$ represents the corresponding quadrature weight. The complete procedure for
matrix-free evaluation is presented in Algorithm~\ref{alg:mf_generic}, where the code for computing the crucial product
$\mathbb{L} :\Grad \corr$ is what we aim to generate automatically.
\begin{algorithm2e}[!ht]
	\SetKwInOut{Input}{Given}
	\SetKwInOut{Output}{Return}
	\Input{  $\currReal$ -- vector representing   $\currFE$ \\
		$ \corrReal $ -- input vector representing  $\corr$ }
	\Output{ $\resultReal_i =  \mathcal{K}(\bar{\bfm{u}}; \corr, \testFunc ) \quad \forall \testFunc \in
			\mathbb{B} $ }
	$\resultReal=0$  \tcp*{zero destination vector }
	\ForEach{ element $K \in \Omega^h$ }{
		gather local vector values on this element\;
		evaluate at	each quadrature point: \\
		\enskip $\Grad \bar {\bfm{u}} ,\; \; \Grad \corr$  \tcp*{Sum factorization}
		\ForEach{quadrature point $q$ on $K$ \tcp*{Quadrature loop}}{
			compute $\bfm{G} = \mathbb{L} :\Grad  \corr $ \tcp*{AceGen-generated}
			queue $\bfm{G}$ for contraction\;
		}
		evaluate queued contractions:
		$\bfm{G} : \Grad\bm{\phi}_i $\tcp*{Sum factorization}
		scatter results to $\resultReal$
	}
	\caption{Matrix-free evaluation of the tangent operator}
	\label{alg:mf_generic}
\end{algorithm2e}

A naive implementation of this procedure typically requires $\mathcal{O}((p+1)^{2d})$ operations for a degree $p$
polynomial basis in $d$-dimensional space. This complexity comes from the fact that evaluating the function
values/gradients at each quadrature point involves looping over all basis functions, leading to quadratic growth in
computational cost
relative to the number of degrees of freedom per element. However, matrix-free methods typically use sum
factorization~\cite{orszag1980spectral,kronbichler2019fast,fischer2020scalability} to bring down the cost of evaluating
the solution
gradients \(
\Grad  \corr\) at quadrature points.
This is achieved by exploiting the tensor-product structure and performing multidimensional evaluation as a series of
1D operations, the complexity is thus reduced to $\mathcal{O}((p+1)^{d+1}) = \mathcal{O}( N_c  \sqrt[d]{N_c} )$,
where
$N_c
	=(p+1)^d$ is the number of degrees of freedom per element.  When	evaluating  gradients, the gradients of
the
mapping
from the reference cell to the physical space are required, as the derivative has to be scaled by the Jacobian matrix
of
the transformation. These Jacobian matrices are precomputed and
stored within the matrix-free data structures.

Sum factorization leads to substantial computational cost savings — for example by a factor of around 16 when using a
moderate polynomial degree of $p=3$ in 3D compared to a full-matrix operation on a single cell with cost
$\mathcal{O}(N_c^2)$. This growth in complexity for a single cell is typically transferred to
the global matrix--vector multiplication if a sparse matrix is used. However, this reduction can be offset
by operations done by the quadrature loop, which, although of linear complexity, can still be expensive. The efficiency
of the quadrature
loop is critical, as a large fraction of the computational time is spent there, as shown in~\cite{davydov2020matrix}.

\subsection{Point-wise evaluation and code generation}
\label{sec:pointwise}

The computational efficiency of the point-wise operations performed at quadrature points
(Algorithm~\ref{alg:mf_generic}) is a critical factor determining the overall matrix-free performance. To handle these
crucial calculations, two primary strategies can be employed, differing fundamentally in their approach to balancing
computational work and memory usage: direct (on-the-fly) computation, where all terms are recalculated as needed, and
partial assembly, which involves pre-computing and caching intermediate results.
Direct computation evaluates all constitutive terms during each operator application; it minimizes memory usage but can
be computationally intensive, especially for complex models.
On the other hand, partial assembly reduces computations during operator application at the cost of increased memory
usage.
Both approaches are broadly applicable, but neither
represents an optimal solution for all scenarios.
The optimal balance depends on the model complexity and hardware characteristics, particularly memory bandwidth.

In~\cite{davydov2020matrix,schussnig2024matrix}, several caching strategies for neo-Hookean elasticity were explored,
it was found that caching
the fourth-order tensor
(exploiting symmetries) was often the most effective strategy for higher polynomial degrees. It was further noted that
this strategy
could
generalize to other models using the spatial tangent tensor, although simpler caching could also be competitive.
In~\cite{schussnig2024matrix}, a similar implementation was tested on newer hardware, finding that
model-specific caching approaches with only scalar quantities provided the best performance.

Let us first detail the direct computation approach, specifically when using the formulation in the reference
configuration. In this case, referring to Eq.~\eqref{eq:DF}, the core term to be computed  is the product of the
tangent stiffness
tensor $\mathbb{L}$ and the referential gradient of the displacement correction $\Delta\bfm{u}$, i.e.,
$\mathbb{L}:\Grad\Delta\bfm{u}$. One can compute $  \mathbb{L} :\Grad  \corr $ on the fly. This requires a minimal
amount of data and is more computationally intense than using pre-computed values. While this might be the natural
choice for simple problems, it may not pay off for more complex ones. Since all the model-dependent data is
recomputed
every time, efficient implementations are crucial for the performance.
We use	AceGen~\cite{korelc2002multi} for finding the expressions of terms at the quadrature points.

\subsubsection{Evaluation on the fly}
\label{sec:pointwise:direct}
In a naive implementation, one would compute $\mathbb{L}$ and double-contract
it with $\Grad\Delta\bfm{u}$. In the context of AD, $\mathbb{L}$ could be obtained by differentiating the strain energy
$\Psi$ twice with respect to the deformation gradient $\bfm{F}$, see Eq.~\eqref{eq:L}. However, thanks to capabilities
of the AD technique
implemented in AceGen, we can use an approach that avoids explicitly forming and contracting
$\mathbb{L}$.
Note that the tangent stiffness tensor $\mathbb{L}$ itself is not needed, only its product with $\Grad\Delta\bfm{u}$,
\begin{equation}
	\label{eq:LG}
	\bfm{G} := \mathbb{L} : \Grad\Delta\bfm{u} ,
\end{equation}
which is to be contracted with the referential gradient of $\delta\bfm{u}$, see Eq.~\eqref{eq:DF}.

% Matrix--vector multiplication can be efficiently implemented in AD by the concept of the \emph{seed matrix}~\cite{griewank2008evaluating}. Let $\bfm{y}(\bfm{x})$ be a set of functions of independent variables $\bfm{x}$.
% Within the AD, the Jacobian is calculated as $\bfm{J}:=\frac{\partial \bfm{y}}{\partial \bfm{x}}=\frac{\partial \bfm{y}}{\partial \bfm{x}}\cdot\bfm{S}$ where $\bfm{S}=\frac{\partial \bfm{x}}{\partial \bfm{x}}=\bfm{I}$ is a seed matrix. 
% By setting the diagonal of the seed matrix to $S_{ii}=v_i$ where $\bfm{v}$ is an arbitrary vector, a matrix--vector multiplication ($\bfm{J}\cdot\bfm{v}$) is computed directly without having to explicitly form the Jacobian matrix and perform the multiplication. The concept of the seed matrix is fundamental to AD and as such can be implemented by most modern AD tools.

%% version updated by JK (and SS) on 20/05
Matrix--vector multiplication can be efficiently implemented in AD by the concept
of the ``seed''~\cite{griewank2008evaluating}. Let $\bfm{y}(\bfm{x})$ be a set of functions of independent variables
$\bfm{x}$. Within the AD, the Jacobian $\bfm{J}:=\frac{\partial \bfm{y}}{\partial \bfm{x}}$ is calculated as
$\bfm{J}=\frac{\partial \bfm{y}}{\partial \bfm{x}} \cdot \frac{\partial \bfm{x}}{\partial \bfm{x}}$, where
$\frac{\partial \bfm{x}}{\partial \bfm{x}}=\bfm{I}$ is a seed matrix. By assuming that $\bfm{x}$ additionally depends
on a fictitious variable $\xi$ with derivatives $\frac{\partial \bfm{x}}{\partial \xi}=\bfm{v}$, the chain rule leads
to $\frac{\partial \bfm{y}}{\partial \xi}=\frac{\partial \bfm{y}}{\partial \bfm{x}} \cdot \frac{\partial
		\bfm{x}}{\partial \xi}=\bfm{J} \cdot \bfm{v}$ where $\bfm{v}$ is an arbitrary seed vector. Thus, using
the AD the
matrix--vector multiplication is performed in a very efficient way without having to explicitly form the Jacobian
matrix and do the multiplication. The concept of the seed is fundamental to AD and as such can be implemented by modern
AD tools.

To arrive at the desired formulation, recall that the first Piola--Kirchhoff stress tensor $\bfm{P}$ depends on the
deformation gradient $\bfm{F}$, thus $\bfm{P}=\bfm{P}(\bfm{F})$, and that $\mathbb{L}$ is defined as the derivative of
$\bfm{P}$ with respect to $\bfm{F}$. Further, assume that the deformation gradient depends on a fictive scalar
variable $\xi$ such that the derivative of $\bfm{F}$ with respect to $\xi$ is equal to $\Grad\Delta\bfm{u}$, so that
\begin{equation}
	\bfm{F} = \bfm{F}(\xi) , \qquad
	\frac{\partial \bfm{F}}{\partial \xi} := \Grad\Delta\bfm{u} .
\end{equation}
Now, using the above assumptions and applying the chain rule, compute the derivative of $\bfm{P}$ with respect
to~$\xi$:
\begin{equation}
	\frac{\partial \bfm{P}}{\partial \xi} \bigg|_{\frac{\partial\bfm{F}}{\partial\xi}=\Grad\Delta\bfm{u}}
	= \frac{\partial \bfm{P}}{\partial \bfm{F}} : \frac{\partial \bfm{F}}{\partial \xi}
	\bigg|_{\frac{\partial\bfm{F}}{\partial\xi}=\Grad\Delta\bfm{u}}
	= \mathbb{L} : \Grad\Delta\bfm{u} ,
\end{equation}
where $\frac{\partial \bfm{F}}{\partial \xi}$ serves as a seed. This formulation yields the desired quantity
$\bfm{G}$, as defined in Eq.~\eqref{eq:LG}.

The above formulation can be implemented in AceGen using the so-called AD exception~\cite{korelc2009automation}.
Assuming that the energy function $\Psi = \Psi(\bfm{H})$ is defined, for gradients $\bfm{H}=\Grad\bar{\bfm{u}}$ and
$\Delta\bfm{H}=\Grad\Delta\bfm{u}$, the evaluation of $\bfm{G}$ is done with the following AceGen/Mathematica code
snippet:
\begin{lstlisting}[
	language=Mathematica, 
	basicstyle=\ttfamily\small\color{black}, 
	mathescape=true,
	backgroundcolor=\color{black!2},
	frame=single,
	rulecolor=\color{black!90!black},
	framerule=1pt,
	frameround=tttt,
	keywordstyle=\color{blue},
	stringstyle=\color{magenta},
	commentstyle=\color{green!60!black}
]
	$\bfm{P} \;  \mathrel{\scriptstyle\bm{\models}}  	 \mathtt{SMSD}[\Psi,\bfm{H}]$;
	$\xi \;\;   \mathrel{\scriptstyle\bm{\models}}  	 	 \mathtt{SMSFictive}[\,]$;
	$\bfm{G}  \,    \mathrel{\scriptstyle\bm{\models}}  	 \mathtt{SMSD}[\bfm{P},\xi,\mathtt{"Dependency"}\to\{\bfm{H},\xi,\Delta\bfm{H}\}]$;
\end{lstlisting}
Here, $\mathtt{SMSD}[\,\cdot\,,\,\cdot\,]$ is a call to the AD procedure that evaluates the
derivative of the first
argument with respect to the second argument, while the option $\mathtt{"Dependency"}$ introduces an AD exception that
intervenes in the AD procedure by overriding the actual dependence existing in the algorithm (here, no dependence, as
$\bfm{P}$ does not depend on $\xi$) by the one specified by this option.

\subsubsection{Partial assembly}
An alternative strategy is to precompute and store $\mathbb{L}$ at quadrature points. This kind of partial assembly
might be beneficial as it does not require any problem-dependent evaluation during the application of the tangent.
However, this
significantly increases the amount of data stored and as a consequence, the algorithm could become memory-bound.

Significant reductions in storage can be achieved by exploiting the symmetries of the fourth-order tensor that are
available when evaluating in the current configuration. Following Eq.~(\ref{eq:DF:cur}), the  application
of the tangent operator to $\corrReal$
%  $\Grad\Delta\bfm{u}$
can be obtained using stored quadrature point data, the fourth-order spatial elasticity tensors $\bbc$
and
the second-order stress tensors $\bm{\sigma}$. We note that, due to the symmetries, $\bbc$ can be stored using only
21 unique real numbers and  $\bm{\sigma}$ requires 6 real numbers in 3D.

The procedure, illustrated in Algorithm~\ref{alg:mf_cached}, is analogous to Algorithm~3 presented
in~\cite{davydov2020matrix}. However, while~\cite{davydov2020matrix} utilized closed-form expressions derived manually
for the quadrature point computations, this work employs automatically generated code via AceGen for the same step
(marked with a comment in the algorithm). The performance of this generated code will be compared against the
hand-written implementation from~\cite{davydov2020matrix}.

\begin{algorithm2e}[!ht]
	\SetKwInOut{Input}{Given}
	\SetKwInOut{Output}{Return}
	\Input{ $ \corrReal $ - input vector representing  $\corr$ \\
		$\bbc$	for each quadrature point\\
		$\bm{\sigma}$ for each quadrature point
	}
	\Output{ $\resultReal_i =  \mathcal{K}(\bar{\bfm{u}}; \corr, \testFunc ) \quad \forall \testFunc \in
			\mathbb{B}$ }
	$\bfm{w}=0$  \tcp*{zero destination vector }
	\ForEach{ element $K \in \Omega^h$ }{
	gather local vector values on this element\;
	evaluate at	each quadrature point: \\
	\enskip $ {\grad}  \corr$  \tcp*{Sum factorization}
	\ForEach{quadrature point $q$ on $K$\tcp*{Quadrature loop}}{
		compute \\
		$\quad \bfm{g}= \bbc : \grads \Delta\bfm{u} + \bm{\sigma} \cdot \grads \Delta\bfm{u}
		$\tcp*{AceGen-generated}
		queue $\bfm{g}$ for contraction\;
	}
	evaluate queued contractions:
	$\bfm{g} : {\grad \phi_i} $\tcp*{Sum factorization}
	scatter results to $\resultReal$
	}
	\caption{
		Matrix-free application of the tangent operator using partial assembly.
		The evaluation at each quadrature point is done using the  cached
		fourth-order spatial elasticity tensor $\bbc$ and Cauchy stress $\bm{\sigma}$,
		cf.\ Algorithm~3 in~\cite{davydov2020matrix}.}
	\label{alg:mf_cached}
\end{algorithm2e}

%%% ===========================================

\section{Performance of matrix-free implementations}
\label{sec:performance}

The ideas discussed above have been implemented in the C++ programming language by
extending the code used by Davydov et al.~\cite{davydov2020matrix}, which builds on the \texttt{deal.II} finite
element 	library~\cite{dealII95}.

We consider a compressible neo-Hookean model and its variant based on the isochoric-volumetric splitting.
For the first model we can directly compare the state-of-the-art hand-written code~\cite{davydov2020matrix} with the
automatically generated code using AceGen. We focus on the comparative performance on	the node level and omit
extensive scalability tests on many nodes, where our algorithm
follows the state-of-the-art, see the results in~\cite{kronbichler2018performance,arndt2020exadg}.
% \mwnote{wording intentional}
{We test various caching strategies that can be employed to find a
balance between computational cost and memory usage towards the goal of a minimal time-to-solution, see
Section~\ref{sec:pointwise} and Section~4 in \cite{davydov2020matrix}.}

\subsection{Setup for performance evaluation}
All the  codes were compiled with \texttt{GCC} 11.4.0 using the following optimization options:
\begin{verbatim}
-O3 -ffast-math -funroll-loops -ffp-contract=fast -march=native
\end{verbatim}
We conduct a series of experiments on a dual-socket AMD EPYC 7282 machine.
We measure the arithmetic floating point processing rate using LIKWID-bench \cite{gruber2022likwid} to describe the
efficiency of the implementation.
The machine's peak performance, obtained via the \texttt{peakflops\_avx} test, is 50 GFLOP/s per core. With a total
of 32 cores, a peak performance of 747 GFLOP/s has been measured.
The machine is equipped with 256 GB of DDR4 memory, with a memory bandwidth of 320 GB/s.
The experiments in parallel are conducted using 32 MPI ranks, which is the maximum
number that can be used without oversubscribing the machine.

\begin{figure}[!ht]
	\centering
	\begin{subfigure}[b]{0.3\textwidth}
		\centering
		\includegraphics[width=\textwidth]{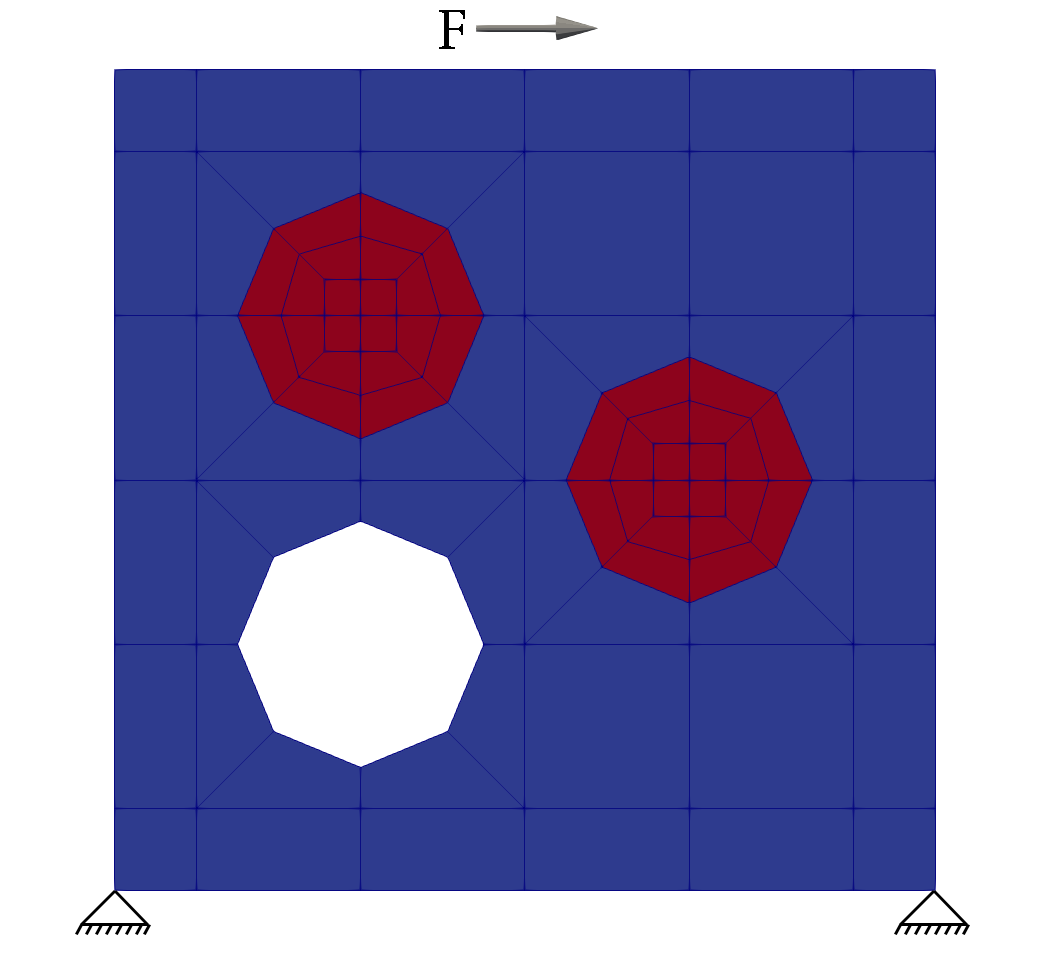}
		\caption{2D coarse mesh}
	\end{subfigure}
	\begin{subfigure}[b]{0.3\textwidth}
		\centering
		\includegraphics[width=\textwidth]{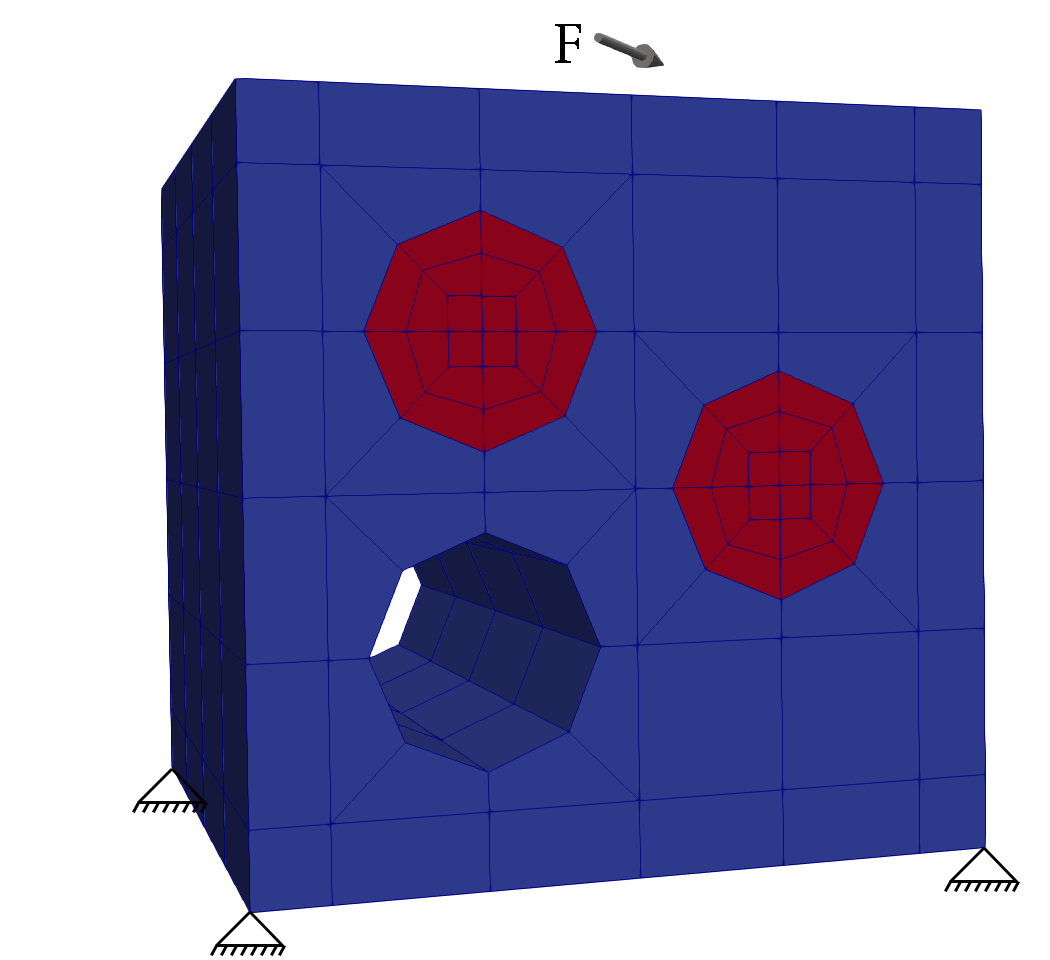}
		\caption{3D coarse mesh}
	\end{subfigure}
	\caption{Discretization of the heterogeneous structure at the coarsest mesh level and the prescribed boundary
		conditions. Both the figure and mesh are taken from the paper by Davydov et al.
		\cite{davydov2020matrix}.}%
	\label{fig:miehe}
\end{figure}

Figure~\ref{fig:miehe} illustrates
the geometry and the discretization of the heterogeneous structure from the test case used in~\cite{davydov2020matrix}
at the coarsest mesh level.	 The 2D structure
consists of a square matrix material (depicted in blue) with a hole and two circular inclusions (depicted in red). The
inclusions are 100 times stiffer than the surrounding material. The 3D geometry is an extrusion of the 2D geometry,
resulting in a cube with edge length of 1000~mm. We use 5 elements in the extrusion direction for the coarse 3D mesh.
The matrix material has a Poisson ratio of 0.3 and a shear modulus of $\mu = 0.4225 \times 10^6 \;
	\text{N}/\text{mm}^2$.

The bottom surface is fixed, while a distributed load is applied incrementally (in 5 loading steps) at the top. The
loading is along the $(1, 0)$ direction and possesses an intensity of $12.5 \times 10^3 \; \text{N}/\text{mm}^2$ for
the
2D problem. For
the
3D problem, the loading is along the $(1, 1, 0)$ direction and possesses an intensity of $12.5 \sqrt{2} \times 10^3 \;
	\text{N}/\text{mm}^2$. The Newton solver uses a displacement tolerance of
$10^{-5}$ and a residual force tolerance of $10^{-8}$, while the linear solver is set to a relative threshold of the
residual of
$10^{-6}$.

\subsection{Model problem: neo-Hookean hyperelasticity}
\label{sec:neo_hookean}

The neo-Hookean model is a widely adopted hyperelasticity model in the study of rubber-like materials due to its
simplicity and good predictive capabilities. In its compressible form, the neo-Hookean elastic strain energy is
expressed as:
\begin{equation}\label{eq:hypref}
	\Psi=\frac{\mu}{2} \left( \tr \bfm{C} - \tr \bfm{I} -2 \log J \right)+ \lambda \log^2 J.
\end{equation}
where $\mu$ is the shear modulus and $\lambda$ is the Lam\'e constant. We refer
to~\cite{ogden1972large,treloar1976mechanics,simo1984remarks,simo1998numerical} for general discussions on the
formulation
of the elastic strain energy for compressible and incompressible hyperelastic materials. The model presented in
Eq.~\eqref{eq:hypref} will be used as the main element for the evaluation of our matrix-free implementation.
This hyperelasticity model was also examined in the studies by~\cite{davydov2020matrix,schussnig2024matrix}, where the
tangent stiffness tensor was derived manually. The simple form of this energy expression enabled several cache
optimizations in these
works. These results provide a baseline for comparison against the automatically generated code approach. Throughout
this text, the model presented in Eq.~\eqref{eq:hypref} will be referred to
as
the \emph{compressible neo-Hookean model}.

For the implementation, we generate code using AceGen.
We note that the evaluation using partial assembly could be further optimized as nearly half of the values of the
stored tensors are zero. While this optimization is easy to implement with the presented tools, we do not
consider it in this work. The reason is that the evaluation presented here using the partial assembly does not depend
on the specific material model, and we want to keep the implementation and the results as general as possible.

\subsubsection{Efficiency of evaluation at quadrature points}
As the first test, we run the program in serial
and measure the execution time of a matrix--vector multiplication with the matrix-free strategy. This setting ensures
that arithmetic costs are most clearly identified, as synchronization overhead is eliminated and memory access costs
through
a shared interface are minimized.
We use $\mathbb{Q}_2$ elements on a relatively coarse mesh (one refinement level) in
3D. Each of the four cores of the CPU	had access to 16 MiB of L3 cache (64 MiB in total). This capacity is sufficient
to accommodate most of the data required by the various caching strategies, ensuring that cache size did not pose a
limiting factor in the performance comparisons.

% \mkwarning{I am not sure 64 MiB are usable -- I believe one core will at most
% 	use the L3	from a core complex/CXX of 16 MiB {https://en.wikipedia.org/wiki/Zen\_2}}

Table~\ref{tab:flops_NH} presents timings, data storage requirements for matrix-free tangent application, and floating
point operations (measured via LIKWID~\cite{gruber2022likwid}).
It compares our automatically generated (AD) strategies with hand-written ones from~\cite{davydov2020matrix}.
The AD strategies include \tensorAD (storing the fourth-order tensor, see Algorithm~\ref{alg:mf_cached}) and \none
(on-the-fly computation, see Algorithm~\ref{alg:mf_generic}).
Among the hand-written strategies from~\cite{davydov2020matrix}: {\tensorS} (Algorithm~3 in~\cite{davydov2020matrix})
is analogous to \tensorAD by also storing the fourth-order tensor.
While our \none approach has no exact hand-written counterpart, {\scalarRef} (see the Appendix
in~\cite{davydov2020matrix}) is the most comparable, as it also involves minimal data storage by caching only $\log(J)$
and performing evaluation in the referential configuration.
For completeness, other hand-written strategies considered are {\scalarCurr} (Algorithm~1 in~\cite{davydov2020matrix}),
which uses the same cached data as \scalarRef but evaluates in the deformed configuration, and the model-specific
{\tensorFAST} (Algorithm~2 in~\cite{davydov2020matrix}), caching a second-order tensor.
% The specific details and correspondences of these strategies are elaborated when discussing the table results below.

\begin{table}
	\caption{Performance metrics for different caching strategies in the quadrature loop  for compressible
		neo-Hookean model.
		Timings of vmult, floating point operations (FLOPs) per point, processing rate in GFLOP/s, and total
		cache size  for Q2 elements in 3D on a mesh
		with one refinement level (75,072 degrees of freedom).
		\label{tab:flops_NH}}
	\begin{centering}
		\input{results/NHDavydov/small_table}
		\par\end{centering}
\end{table}

The first three rows in Table~\ref{tab:flops_NH} present the results for the automatically generated code.
In the first row, we show the results obtained using a naive approach to AD, where the tangent stiffness tensor is
explicitly formed through standard differentiation and then contracted with the gradient, rather than using the
efficient AD exception technique described in Section~\ref{sec:pointwise:direct}. The second row (\none) implements
the computation of the product $\mathbb{L}:\Grad\Delta\bfm{u}$ directly
without explicitly forming $\mathbb{L}$, resulting in significantly improved performance. The third row shows
results for storing the fourth-order tensor (\tensorAD).
We observe that even for such a small problem, strategy \tensorAD performs
only slightly better than the one without caching (\none), despite its lower number of operations per quadrature point.
We expect that this effect can be attributed to the significantly larger cache size.

The next four rows in Table~\ref{tab:flops_NH} detail the performance of the hand-written strategies from
\cite{davydov2020matrix}.
The {\scalarRef} strategy, despite having fewer operations per quadrature point than the AD-based \none approach, is
slower.
The {\scalarCurr} strategy also shows lower performance than the AD-based strategies.
For both {\scalarRef} and {\scalarCurr}, this difference can be attributed to a lower processing rate in the
hand-written codes, likely due to compiler-generated overhead and more computationally expensive operations.
The processing rate of the AD-generated \none strategy is nearly twice that of the hand-written \scalarCurr strategy.

The hand-written {\tensorS} strategy is clearly slower and requires significantly more data than its AD counterpart
(\tensorAD), as not all tensor symmetries are exploited in the manual implementation. The automatically generated
version also achieves an over 50\% higher processing rate. (For reference, the sum factorization part of the
matrix--vector product was performed at around 29.8 GFLOPS/s, although for very high order elements, the peak
performance
was around 44 GFLOPS/s.)

Finally, {\tensorFAST}, the fastest hand-written approach, has the lowest
number of operations per quadrature point and a cache size only twice than that of {\scalarRef}. This approach is not
generalizable as it relies on the specific form of the compressible neo-Hookean model,
Eq.~\eqref{eq:hypref}, leading to a particularly simple spatial tangent tensor, see Eqs.~(15) and (16)
in~\cite{davydov2020matrix}.

\subsubsection{Parallel performance}
For a more realistic balance between compute capability and memory bandwidth, we here test our program in parallel by
running the code with 32 MPI processes across a range of polynomial
degrees and refinement levels. To ensure a balanced comparison, the number of degrees of freedom is maintained within
the
same order of magnitude, as detailed in Tables \ref{tab:input_parameters_2d} and \ref{tab:input_parameters_3d}. For
instance, in 3D, when using polynomial degree $p=1$, we apply four refinement levels, while for $p=4$, we apply only
one
refinement level to maintain a similar computational scale.

\begin{table}[!htb]
	\caption{ {Parameters for the benchmark: $p$ is the polynomial degree,
				$q$ is the number of quadrature points in 1D, $N_{\text{gref}}$ is the number of global
				mesh refinements,
				$N_{el}$ is the number of elements and $N_{\text{DoF}}$ is the number of DoFs.}}
	\begin{subtable}{.49\linewidth}
		\caption{2D}
		\label{tab:input_parameters_2d}
		\centering
		\begin{tabular}{ccccc}
			\hline
			$p$ & $q$ & $N_{\text{gref}}$ & $N_{el}$  & $N_{\text{DoF}}$ \\
			\hline
			1   & 2   & 7                 & 1,441,792 & 2,887,680        \\
			2   & 3   & 6                 & 360,448   & 2,887,680        \\
			3   & 4   & 5                 & 90,112    & 1,625,088        \\
			4   & 5   & 5                 & 90,112    & 2,887,680        \\
			5   & 6   & 5                 & 90,112    & 4,510,720        \\
			6   & 7   & 4                 & 22,528    & 1,625,088        \\
			7   & 8   & 4                 & 22,528    & 2,211,328        \\
			8   & 9   & 4                 & 22,528    & 2,887,680        \\
			\hline
		\end{tabular}
	\end{subtable}
	\begin{subtable}{.49\linewidth}
		\caption{3D}
		\label{tab:input_parameters_3d}
		\centering
		\begin{tabular}{ccccc}
			\hline
			$p$ & $q$ & $N_{\text{gref}}$ & $N_{el}$  & $N_{\text{DoF}}$ \\
			\hline
			1   & 2   & 4                 & 1,441,792 & 4,442,880        \\
			2   & 3   & 3                 & 180,224   & 4,442,880        \\
			3   & 4   & 2                 & 22,528    & 1,891,008        \\
			4   & 5   & 2                 & 22,528    & 4,442,880        \\
			\hline
		\end{tabular}
	\end{subtable}
\end{table}

\begin{figure}
	\begin{centering}
		\begin{tabular}{cc}
			\begin{tikzpicture}
				\node[anchor=south west,inner sep=0] (image)
				{\includegraphics[scale=\plotscale]{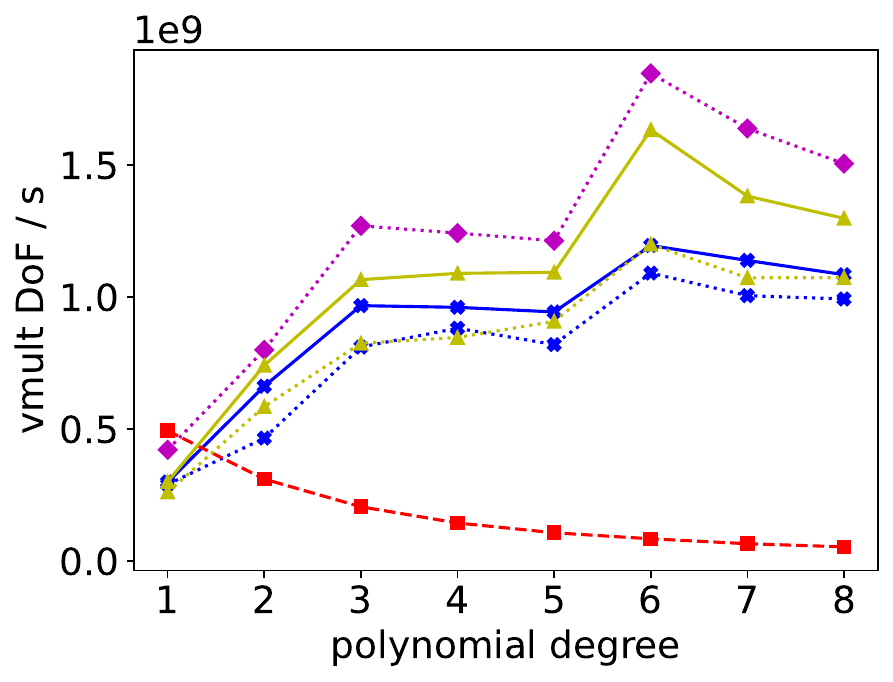}};
				\node [anchor=north east, color=black, inner sep=1pt,
					xshift=-15.3em, yshift=-2em] at (image.north east) {\boxed{{\rm 2D}}};
			\end{tikzpicture} &
			\begin{tikzpicture}
				\node[anchor=south west,inner sep=0] (image)
				{\includegraphics[scale=\plotscale]{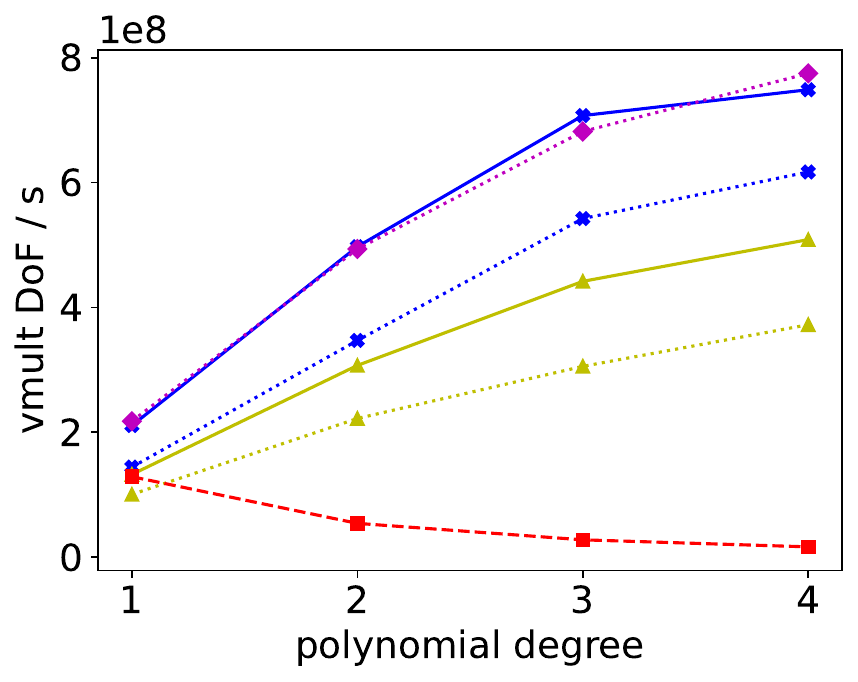}};
				\node [anchor=north east, color=black, inner sep=1pt,
					xshift=-15.3em, yshift=-2em] at (image.north east) {\boxed{{\rm 3D}}};
			\end{tikzpicture} \\

			\multicolumn{2}{c}{\includegraphics[scale=0.5]{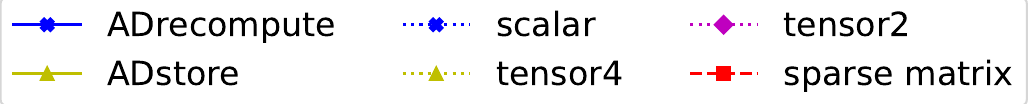}}
		\end{tabular}

		\caption{Measured throughput of matrix--vector operator evaluation for the compressible neo-Hookean
			model.
			The processing rate is expressed in DoFs/second. The data is shown for 2D (left) and 3D
			(right). The results obtained with automatically generated code are depicted with solid lines,
			while the ones obtained with hand-written code~\cite{davydov2020matrix} are shown with dotted
			lines. The sparse-matrix vmult is shown by red dashed line.
		}
		\label{fig:NHDavydovTiming}
	\end{centering}
\end{figure}

Below we compare the timing of a single matrix--vector multiplication for various evaluation methods, including the
``classical''
way, i.e., involving a sparse matrix operations from the Epetra
package of the Trilinos project~\cite{heroux2005overview}.  As a
normalized measure of performance, we record
the processing rate in DoFs/second. The obtained results are depicted in
Figure~\ref{fig:NHDavydovTiming}.
We have tested all the formulations listed in Table~\ref{tab:flops_NH}, however, here we skip the results for
\scalarRef as its performance is comparable to, but slightly worse
than, that of \scalarCurr. The results obtained with automatically generated code are depicted with solid lines,
while the ones obtained with hand-written code are shown with dotted lines.
We observe that for $p>1$ all matrix-free operators outperform the sparse matrix one, with the gap growing with the
degree $p$.
This can be attributed to the more favorable ratio between the number of degrees
of freedom and the number of quadrature points per cell at higher $p$.
This can be accounted to exploiting the
tensor-product structure of the shape functions, as discussed in Section~\ref{sec:mf_workflow}. For degree $p=4$ in 3D,
matrix-free evaluations are up to 40 times faster than the sparse-matrix approach. When using linear
shape functions the sparse matrix approach is still slower than the matrix-free one by a factor of 4 in 3D.

In 2D the hand-written approach with caching a second-order tensor (\tensorFAST) is the best, closely followed by
automatically generated
code that caches the fourth-order tensor (\tensorAD). The superior performance of caching
approaches is expected, as in 2D the data size is smaller than in 3D.
In contrast, in 3D, the operator without caching performs
comparably to the hand-written implementation that caches a second-order tensor. We recall
that the latter implementation is limited to the considered model.

\begin{comment}
To explain the increase of the efficiency with the polynomial degree we first present the breakdown of the
matrix--vector multiplication in \mwerror{Plot the breakdown}{Figure???}. \mknote*{Can we say this? A CPU will execute
	work possibly out of order, pre-fetching of data access can overlap some regions. I think it would make sense
	to
	present the number of instructions spent at quadrature points versus the rest. In any case, I am curious.}{It
	can be
	observed that the majority of time
	is spent on the quadrature loop even for higher polynomial degrees.} Hence, we attribute that gain to the
decreasing
ratio of the number of total quadrature points to degrees of freedom for increasing $p$.
\end{comment}

\begin{comment}
To explain the increase of the efficiency with the polynomial degree we first present the breakdown of the
matrix--vector multiplication in \mwerror{Plot the breakdown}{Figure???}. \mknote*{Can we say this? A CPU will execute
	work possibly out of order, pre-fetching of data access can overlap some regions. I think it would make sense
	to
	present the number of instructions spent at quadrature points versus the rest. In any case, I am curious.}{It
	can be
	observed that the majority of time
	is spent on the quadrature loop even for higher polynomial degrees.} Hence, we attribute that gain to the
decreasing
ratio of the number of total quadrature points to degrees of freedom for increasing $p$.
\end{comment}

To explain the increasing efficiency for higher polynomial degrees let us consider a Cartesian mesh with $n_{c}$ cells
in
each direction. For a finite
element of degree $p$, we have $n_{\text{dof}}=(p n_c+1)^d$ degrees of freedom and { $ \left(n_{c} (p+1)\right)^d$
		quadrature points}. Then the ratio between the total number of quadrature points to degrees of freedom
is
$$ \Bigg( \frac{(p+1)n_c}{p n_c +1} \Bigg)^d \approx \Bigg( \frac{p+1}{ p } \Bigg)^d.
$$
This ratio for $p=1$ is  4 in 2D and 8 in 3D, and approaches 1 for higher degrees $p$, see
also~\cite{kronbichler2012generic}.

\begin{figure}
	\begin{centering}
		\begin{tabular}{cc}
			\begin{tikzpicture}
				\node[anchor=south west,inner sep=0] (image)
				{\includegraphics[scale=\plotscale]{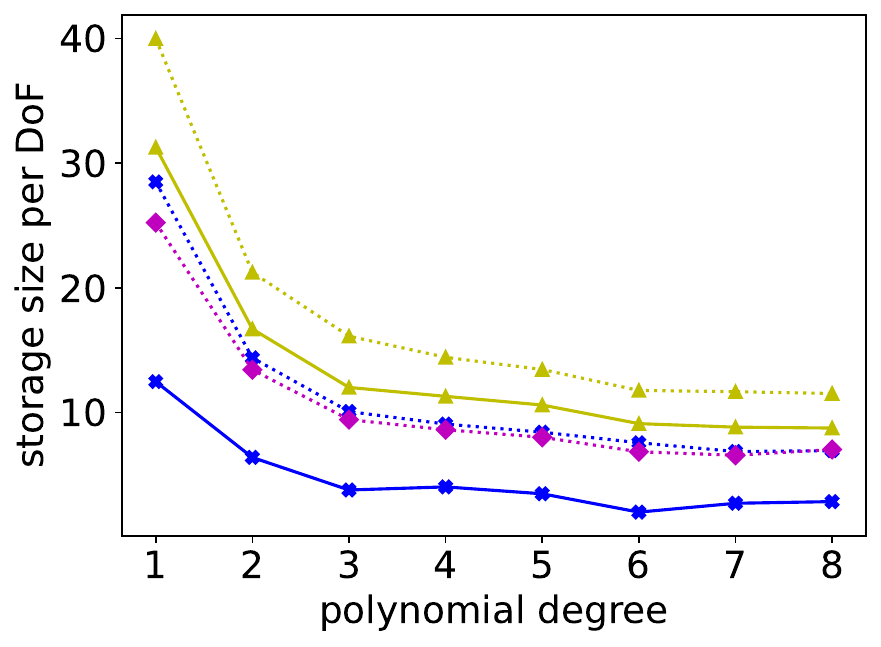}};
				\node [anchor=north east, color=black, inner sep=1pt,
					xshift=-1em, yshift=-1em] at (image.north east) {\boxed{{\rm 2D}}};
			\end{tikzpicture} &
			\begin{tikzpicture}
				\node[anchor=south west,inner sep=0] (image)
				{\includegraphics[scale=\plotscale]{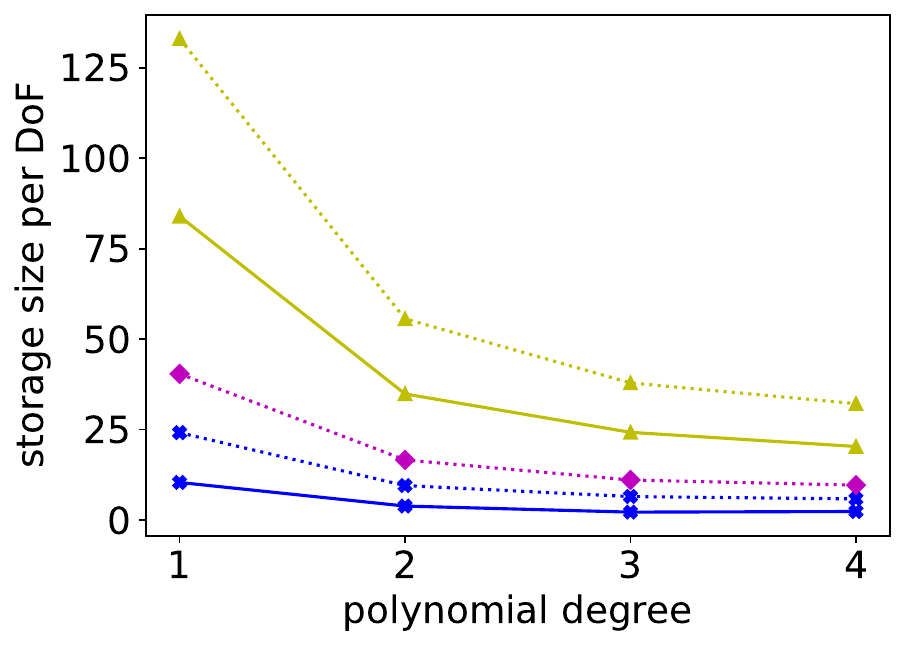}};
				\node [anchor=north east, color=black, inner sep=1pt,
					xshift=-1em, yshift=-1em] at (image.north east) {\boxed{{\rm 3D}}};
			\end{tikzpicture} \\

			\multicolumn{2}{c}{\includegraphics[scale=0.5]{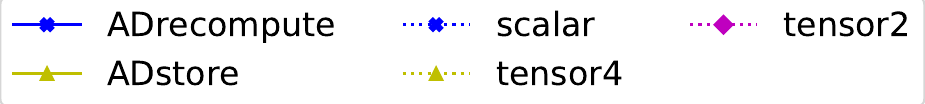}}
		\end{tabular}
		\caption{Memory requirements per degree of freedom for matrix--vector operator application for
			the compressible neo-Hookean
			model. The storage size is expressed in the number of floating point numbers  per DoF. The data
			is
			shown for 2D (left) and 3D (right). The results obtained with automatically generated code are
			depicted with solid lines, while the ones obtained with hand-written
			code~\cite{davydov2020matrix} are shown with dotted
			lines.}
		\label{fig:NHDavydovMem}
	\end{centering}
\end{figure}

In Figure~\ref{fig:NHDavydovMem} we plot the storage size required for the	application
of the	matrix--vector operator. We express the size in the number of floating point numbers per degree of freedom.
Especially in 3D, it is visible that memory usage impacts computing time as the storage of
fourth-order tensors is the slowest strategy, even though it involves the lowest number of computations.
The \tensorAD caching strategy requires storing 27 numbers per quadrature point, in addition to the standard
quadrature-point data needed for any matrix-free operator: Jacobian matrices of the transformation (9 numbers) leading
to total storage of 36 numbers. Note that there are also 3 degrees of freedom per node,
meaning that the minimal storage per degree of freedom is 36/3=12, while the observed ratio for \tensorAD with $p=4$
in 3D is 20.3.	As a reference, at degree $p=4$, the sparse-matrix approach requires 51 times as much memory, and over
750
times more memory than the matrix-free operator without intermediate result storage.

\begin{figure}
	\begin{centering}
		\begin{tabular}{cc}
			\begin{tikzpicture}
				\node[anchor=south west,inner sep=0] (image)
				{\includegraphics[scale=\plotscale]{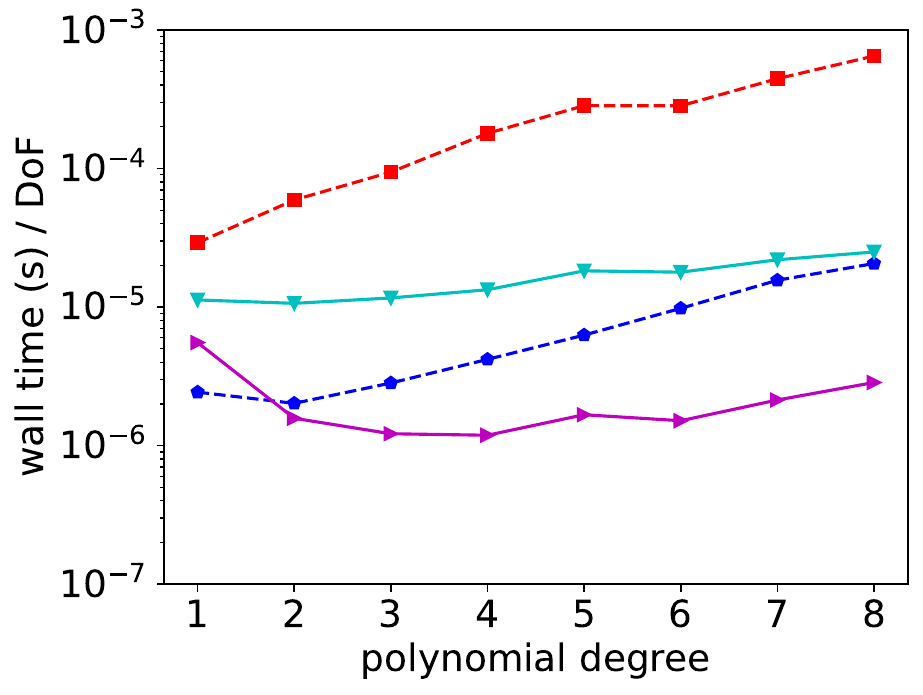}};
				\node [anchor=south east, color=black, inner sep=1pt,
					xshift=-1em, yshift=3em] at (image.south east)
				{\boxed{{\rm 2D}}};
			\end{tikzpicture} &
			\begin{tikzpicture}
				\node[anchor=south west,inner sep=0] (image)
				{\includegraphics[scale=\plotscale]{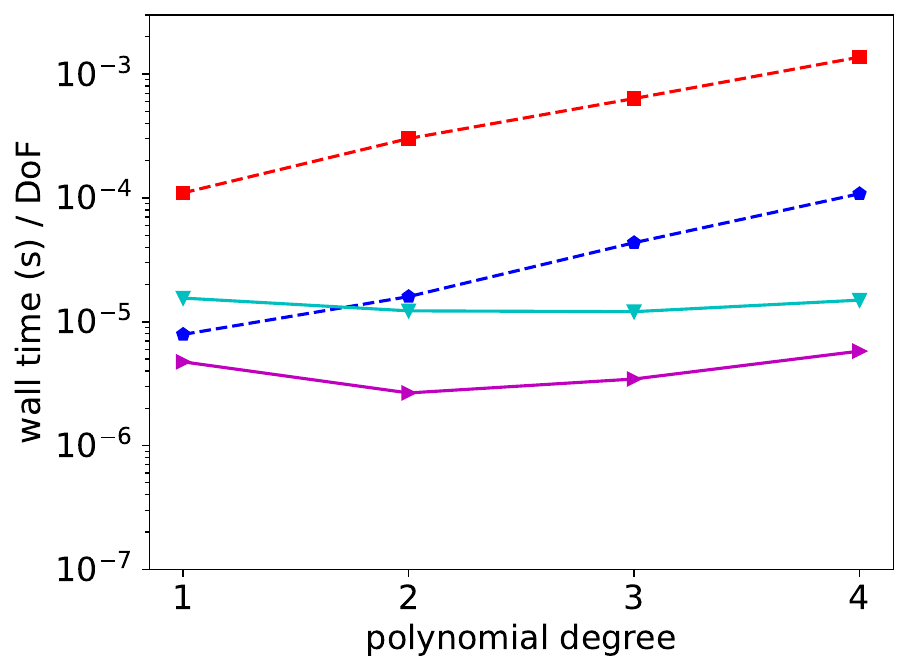}};
				\node [anchor=south east, color=black, inner sep=1pt,
					xshift=-1em, yshift=3em] at (image.south east)
				{\boxed{{\rm 3D}}};
			\end{tikzpicture}
			\\

			\multicolumn{2}{c}{\includegraphics[scale=0.55]{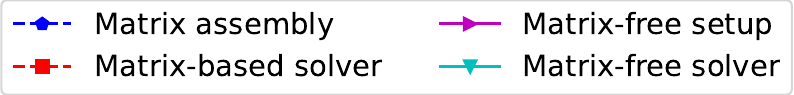}}
		\end{tabular}
		\caption{Comparison of time to solution for matrix-free and sparse matrix approaches across different
			polynomial degrees in 2D and 3D for the compressible neo-Hookean model.
			Computations using the
			sparse-matrix approach are shown with dashed lines, while the matrix-free
			approach is shown with
			solid lines.
		}
		\label{fig:time_to_solution}
	\end{centering}
\end{figure}

Finally, we also plot the time to solution for matrix-free versus sparse matrix approaches. The results,
depicted in Figure~\ref{fig:time_to_solution}, clearly demonstrate the superior performance of matrix-free methods. The
matrix-free approach consistently outperforms the sparse-matrix-based one, with the gap widening as the polynomial
degree increases. For polynomial degree $p=4$ in 3D, the matrix-free solver is 80 times faster than the
sparse-matrix approach.

We observe that the time to solution does not follow the same trend as the processing rate
for the matrix-free operator. This is due to the dependence of the preconditioner on the polynomial degree, which
leads to an increase in the number of iterations for the solver for higher degrees. This issue is associated with the
smoother and can be resolved by using a more sophisticated
smoother~\cite{wichrowski2024smoothers}.%76 vs 6170 s for 15 total solves

%76 vs 6170 s for 15 total solves

\subsection{Neo-Hookean model with isochoric-volumetric split}

Another popular variant of the neo-Hookean hyperelasticity model is the one that splits the elastic strain energy into
isochoric and volumetric parts.
% \ssnote{do we need this sentence? or any reference supporting the statement?}\redtext{This split is particularly
% 	advantageous for the computer implementation of nearly-incompressible materials, as it allows to more
% 	effectively
% 	handle the nearly-incompressible deformation behavior.}
Among the various formulations, the following form of the elastic strain energy is adopted here
\cite{simo1998numerical}
\begin{equation}\label{eq:hypcomplex}
	\Psi=\frac{\mu}{2}(\tr \bar{\bfm{C}}- \tr \bfm{I}) + \frac{\kappa}{2} \left( \frac{1}{2} (J^2-1) - \log J
	\right),
\end{equation}
where $\bar{\bfm{C}}=J^{-\frac{2}{3}}\bfm{C}$ denotes the isochoric part of the right Cauchy--Green tensor $\bfm{C}$
and
$\kappa$ is the bulk modulus. This model is henceforth referred to as \emph{split neo-Hookean}.

Our rationale for selecting this model stems from the additional complexities it introduces in the
evaluation of the residual vector and tangent matrix. This model was also considered in recent
work~\cite{schussnig2024matrix} on matrix-free elasticity solvers, where an explicit derivation of the tangent operator
was provided.

\subsubsection{Performance of matrix-free implementation}
With these preliminary observations in mind, our
goal is to investigate how the inherent complexities of the model translate to the matrix-free implementation. We first
measure floating point operations per
quadrature point and their execution rate as for the previous model. {Table~\ref{tab:flops_NHSplit} summarizes these
		results.}

\begin{table}
	\begin{centering}
		\caption{Performance metrics for different caching strategies in the quadrature
			loop for the split neo-Hookean model.
			Timings of vmult, floating point operations (FLOPs) per point, processing rate in
			GFLOP/s, and total cache size for Q2 elements in 3D on a mesh with one
			refinement level (75,072 degrees of freedom).
			\label{tab:flops_NHSplit}}
		\begin{tabular}{ l |  c  c  c  c }
			\hline
			% &            & \tabularnewline
			Formulation    & Timing [ms] & FLOPs per point & {[}GFLOP/s{]} & Cache
			{[}Mb{]}    \tabularnewline
			\hline
			\noneSplit     & 3.91        & 3634            & 27.0          & 2 \tabularnewline
			\tensorADsplit & 3.27        & 1566            & 19.5          & 18
			\tabularnewline
			\hline
		\end{tabular}
		\par\end{centering}
\end{table}

\begin{figure}
	\begin{centering}
		\begin{tabular}{cc}
			\begin{tikzpicture}
				\node[anchor=south west,inner sep=0] (image)
				{\includegraphics[scale=\plotscale]{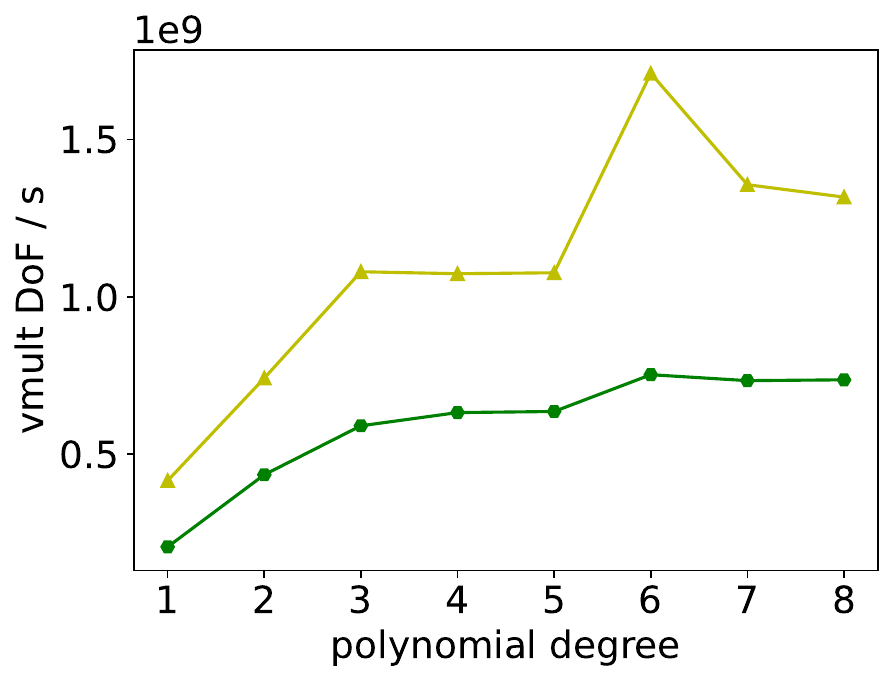}};
				\node [anchor=north east, color=black, inner sep=1pt,
					xshift=-1em, yshift=-11em] at (image.north east) {\boxed{\rm 2D}};
			\end{tikzpicture} &
			\begin{tikzpicture}
				\node[anchor=south west,inner sep=0] (image)
				{\includegraphics[scale=\plotscale]{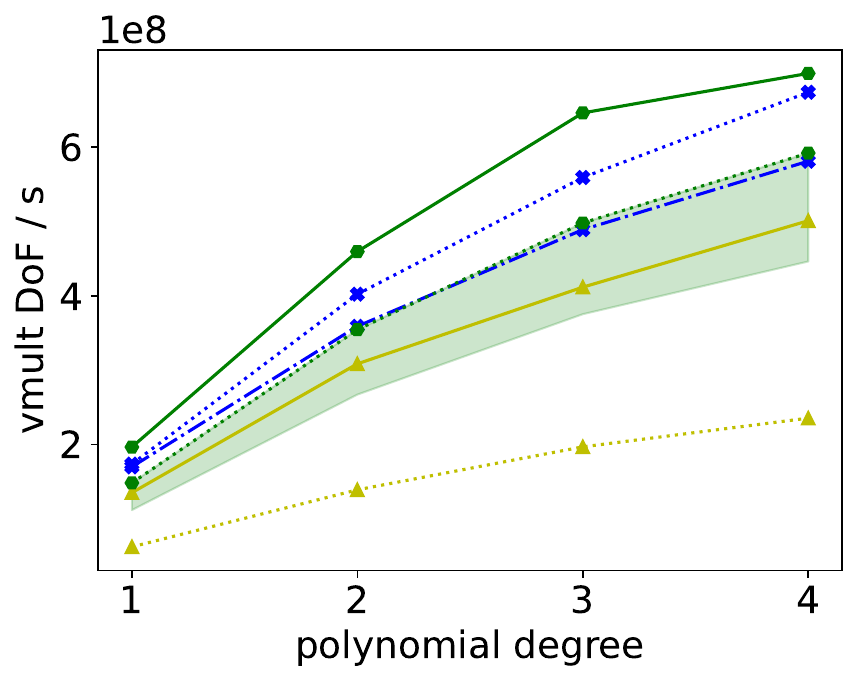}};
				\node [anchor=north east, color=black, inner sep=1pt,
					xshift=-1em, yshift=-11em] at (image.north east) {\boxed{\rm 3D}};
			\end{tikzpicture} \\

			\multicolumn{2}{c}{\includegraphics[scale=0.5]{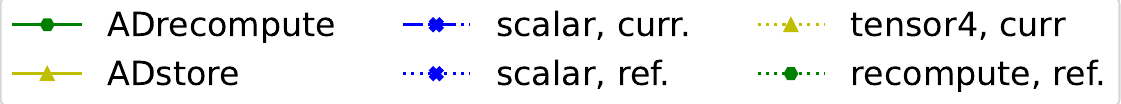}}
		\end{tabular}
		\caption{
			Measured throughput  during application of matrix--vector operator for the split neo-Hookean
			model.
			The processing rate is expressed in DoFs/second. The data is shown for 2D (left) and 3D
			(right). The results obtained with automatically generated code are depicted with solid lines,
			while the estimates for the ones obtained with hand-written code~\cite{schussnig2024matrix} are
			shown with dotted lines.
			To show the uncertainty, we indicate possible variations of \emph{recompute all} throughput
			with the green area.}
		\label{fig:NHSplitTiming}
	\end{centering}
\end{figure}

Following our previous approach, we evaluate the solver performance across various grid sizes and polynomial degrees by
measuring the times of matrix--vector operations. The corresponding results are shown in
Figure~\ref{fig:NHSplitTiming}. Since the
evaluation
timing and memory usage for the \tensorAD caching strategy are model-independent, we use these as reference
values. The results confirm our earlier findings: caching with \tensorAD is advantageous in 2D, while
on-the-fly computation proves more efficient in 3D. This of course will change for more complex models, where
storing the fourth-order tensor is advantageous, as it decouples the tangent evaluation from the specific model.

For this model, a direct performance comparison between our automatically-generated code and hand-written
implementations is not feasible. Instead, we use our implementation of the compressible neo-Hookean model as a
reference and rescale the throughput reported in~\cite{schussnig2024matrix} to account for differences in hardware and
test problems. Specifically, we compute the ratio between our measured throughput and those reported
in~\cite{schussnig2024matrix} for both scalar caching and fourth-order tensor caching scenarios. To ensure a
conservative estimate in each case, we take the maximum ratio across all polynomial degree $p$ considered. For the
\emph{recompute all} strategy reported in their work, we adopt the conservative ratio resulting from scalar caching.
Note that this approach likely results in an overestimation of the implementation performance
for~\cite{schussnig2024matrix}. To reflect the uncertainty and the range of performance variation due to different
polynomial degrees, we depict the resulting spread in Figure~\ref{fig:NHSplitTiming} as a green area, indicating the
range of possible performance ratios. Nevertheless,
our implementation of {\none} still outperforms even these generous estimates for the hand-written code.

\FloatBarrier
\section{Conclusion and outlooks}

This work demonstrates that automating code generation for matrix-free methods in finite-strain elasticity, through
tools like AceGen, significantly enhances development efficiency and yields computationally superior code. By
leveraging automatic differentiation combined with symbolic and stochastic simplifications, the generated code not only
surpasses CPU processing speeds of traditional hand-written codes but is also inherently less prone to human error,
leading to more robust, efficient, and maintainable implementations.

Direct comparisons for the neo-Hookean model reveal that the automatically generated code outperforms its hand-written
counterpart. Our automated approach also demonstrates excellent performance for the more complex neo-Hookean model with
an isochoric-volumetric split, highlighting its versatility.

Our investigation into caching strategies shows that, while on-the-fly computation with minimal data storage is optimal
for the 3D models studied, storing intermediate results is beneficial in 2D due to smaller tensor sizes. For more
complex constitutive models, computing everything on-the-fly may not be optimal. In such scenarios, caching the
fourth-order tensor emerges as a robust strategy, offering a good compute--storage balance by decoupling tangent
evaluation from the
specific material model and still delivering substantially better performance than sparse-matrix methods.
Looking ahead, memory bandwidth may become a bottleneck on future hardware, a challenge potentially mitigated by
reorganizing
computational operations~\cite{wichrowski2024smoothers}.

Further research avenues include addressing current solver limitations. The reliance of multigrid methods on a
multilevel mesh hierarchy can be restrictive for complex geometries; unfitted approaches like cutFEM, were recently
shown to be
compatible with matrix-free techniques~\cite{hansbo2017cut, bergbauer2024high}, offering a promising alternative.
Moreover, specialized techniques are crucial for applying multigrid to nearly-incompressible solids. This challenge can
be tackled by developing more advanced smoothers or by employing mixed formulations with robust block
solvers~\cite{mardal2011preconditioning}.
These areas remain open for investigation with automatically generated matrix-free operators.

\section*{Acknowledgments}
The authors declare the use of language models (ChatGPT, Gemini, and Claude) to improve the clarity and readability of
the manuscript. All scientific content and technical claims are solely the responsibility of the authors.
MRH and SS acknowledge support from the EU through the EffectFact project (No.\ 101008140) funded within the
H2020-MSCA-RISE-2020 programme, and wish to thank Dr.\ Toma\v{z} \v{S}u\v{s}tar for useful discussions and kind
hospitality while visiting C3M, Ljubljana, Slovenia. MK acknowledges support by the EU through the EuroHPC joint
undertaking Centre of Excellence dealii-X (No.\ 101172493).

\bibliographystyle{ieeetr}
\bibliography{literature}

\end{document}

%% file: preamble.tex
\usepackage{amstext}
\usepackage{stmaryrd}
\usepackage{algpseudocode}
\usepackage[algo2e,linesnumbered,vlined,commentsnumbered,ruled]{algorithm2e}

\usepackage[normalem]{ulem}
\usepackage{graphicx}

\usepackage{array}
\usepackage{verbatim}
\usepackage{booktabs}
\usepackage{calc}
\usepackage{textcomp}
\usepackage{multirow}

\usepackage[all]{xy}

\usepackage{cancel}
\usepackage{mathtools}
\usepackage{placeins}
\usepackage{xspace}

\usepackage{tikz,tikzscale}
\usetikzlibrary{calc}
\usetikzlibrary{shadings}
\usepackage[mode=multiuser,status=draft]{fixme} % Provide note, warning, error, fatal.
\fxsetup{envlayout=color}
\fxsetup{innerlayout=inline}
\fxsetup{targetlayout=colorcb}
\FXProvidesTargetLayout{color}

\fxusetheme{color}
\FXRegisterAuthor{ss}{envss}{SS}
\FXRegisterAuthor{mr}{envmr}{MR}
\FXRegisterAuthor{mw}{envmw}{MW}
\FXRegisterAuthor{mk}{envmk}{MK}

\usepackage{subcaption}
\usepackage[format=plain,indention=.5cm]{caption}

\usepackage{tcolorbox}
\usepackage{listings}
\tcbuselibrary{listings}

\definecolor{mathematicablue}{RGB}{0,0,150}
\lstdefinestyle{mathematica}{
	language=Mathematica,
	basicstyle=\ttfamily\small,
	keywordstyle=\color{mathematicablue}\bfseries,
	commentstyle=\color{gray}\itshape,
	stringstyle=\color{purple},
	showstringspaces=false,
	breaklines=true
}

\def\mesh{\mathbb M}

\newcommand{\Dd}{\text{d}}
\newcommand{\DD}{\text{D}}

\usepackage{bm,upgreek}

\newcommand{\bfm}[1]{{\rm\bf #1}}

\DeclareMathOperator{\tr}{tr}
\newcommand{\fem}[1]{\mbox{\bf \texttt{#1}}}

\newcommand{\corr}{  \,\Delta \bfm{u}\, }
\newcommand{\testFunc}{\bm{\phi}_i}

\newcommand{\currFE}{\, \bar{ \bfm{u}}\, }

\newcommand{\currReal}{\bar{ \fem{U}}}
\newcommand{\corrReal}{\Delta \fem{U}}
\newcommand{\resultReal}{\fem{W}}

\newcommand{\scalarRef}{\texttt{scalar\_ref}\xspace}
\newcommand{\tensorS}{\texttt{tensor4}\xspace}

\newcommand{\none}{\texttt{ADrecompute}\xspace}
\newcommand{\tensorAD}{\texttt{ADstore}\xspace}
\newcommand{\tensorFAST}{\texttt{tensor2}\xspace}
\newcommand{\scalarCurr}{\texttt{scalar}\xspace}

\newcommand{\crap}{\texttt{naive\_AD}\xspace}

\newcommand{\noneSplit}{\texttt{ADrecompute}\xspace}
\newcommand{\tensorADsplit}{\texttt{ADstore}\xspace}

\DeclareMathOperator{\Grad}{Grad}
\DeclareMathOperator{\grad}{grad}
\newcommand{\grads}{\grad^{\rm s}\!}
\usepackage{graphicx} % Required for inserting images

\usepackage{bbm}
\newcommand{\bbc}{\mathbbm{c}}

\usepackage{amsmath} % Good practice for math
\usepackage{amssymb} % Provides \models
\usepackage{bm}      % Provides \bm for better bold math
\usepackage{graphicx} % Provides \scalebox

\usepackage{graphicx} % Required for inserting images

%% file: results/NHDavydov/small_table.tex
\begin{tabular}{ l | c	c c c }
	\hline
	Formulation                           & Timing [ms] & FLOP/point & {[}GFLOP/s{]} & Cache {[}Mb{]}
	\tabularnewline
	\hline
	%%%
	% \ifdefined \GAMM
	% \else
	\crap                                 & 5.56         & 6913       & 33.9          & 2 \tabularnewline
	% \fi
	%%
	\none                                 & 3.43         & 3287       & 25.3          & 2 \tabularnewline
	\tensorAD                             & 3.33         & 1566       & 18.5          & 18 \tabularnewline
	\hline
	\scalarRef \cite{davydov2020matrix}   & 4.71         & 3148       & 14.8          & 4
	\tabularnewline
	%%%
	\ifdefined \GAMM
	\else
	\scalarCurr  \cite{davydov2020matrix} & 4.11         & 2138       & 13.2          & 5 \tabularnewline
	\fi
	%%%
	\tensorS  \cite{davydov2020matrix}    & 4.44         & 1588       & 13.5          & 29 \tabularnewline
	\tensorFAST  \cite{davydov2020matrix} & 2.41         & 1146       & 12.4          & 9 \tabularnewline
	\hline
\end{tabular}